\documentclass[11pt]{article}

\RequirePackage[italian,english]{babel}
\usepackage{authblk}

\title{The Hopf-Lax formula for multiobjective costs with non-constant discount 
  via set optimization%
  }
  
\author{D.~Visetti\footnote{Faculty of Economics and Management, Free University of Bozen-Bolzano, Italy.}}

\RequirePackage[italian,english]{babel} 
\RequirePackage{graphicx}
\RequirePackage{amsmath}
\RequirePackage{amssymb}
\RequirePackage{amsthm}

\usepackage[colorlinks=true, linkcolor=blue, citecolor=blue, urlcolor=blue, filecolor=blue]{hyperref}

\newcommand{\cl}{{\rm cl \,}}

\newcommand{\bs}{\backslash}
\renewcommand{\P}{\ensuremath{\mathcal{P}}}
\newcommand{\A}{\ensuremath{\mathcal{A}}}
\newcommand{\F}{\ensuremath{\mathcal{F}}}

\begin{document}
\selectlanguage{english}

\newcommand{\NN}{{\mathbb N}}
\newcommand{\RR}{{\mathbb R}}
\newcommand{\ZZ}{{\mathbb Z}}
\newcommand{\QQ}{{\mathbb Q}}
\newcommand{\CC}{{\mathbb C}}

\newtheorem{theorem}{Theorem}[section]
\newtheorem{proposition}[theorem]{Proposition}
\newtheorem{lemma}[theorem]{Lemma}
\newtheorem{corollary}[theorem]{Corollary}
\newtheorem{remark}[theorem]{Remark}
\newtheorem{example}[theorem]{Example}
\newtheorem{definition}[theorem]{Definition}

\maketitle

\begin{abstract}
  The minimization of a multiobjective Lagrangian with non-constant 
  discount is studied. The problem is embedded into a set-valued framework 
  and a corresponding definition of the value function is given.   Bellman's optimality 
  principle and Hopf-Lax formula are derived.  The value function is shown 
  to be a solution of a set-valued Hamilton-Jacobi equation.
\end{abstract}


\section{Introduction}

In this paper the following optimization problem is considered
$$
\inf_{y\in A(t,x)}\left(\int_t^Td_t(s)L(\dot y(s))ds+g(y(T))\right)\hspace{6cm} (P)
$$
where $L$ is a vector-valued Lagrangian or utility function, $d_t(s)$ is a discount 
factor and $A(t,x)$ is a set of admissible arcs that take the value $x$ at time $t$.   
Often it is the case that there are more than one function to minimize (for 
example production cost and holding inventory cost).  Sometimes these functions 
can contradict each other, in the sense that trying to minimize one of them leads to 
the increase of another one. It can be sensible to discount the expenditures: 
for example we can think of a discount at a continuous rate $r$ 
$d_t(s)=e^{-r(s-t)}$ or $d_t(s)=e^{-\int_t^s\rho(s)ds}$, where $\rho$ is in 
$L^\infty([0,T])$.   In the Preface of \cite{Caputo} the author writes that 
``Forward-looking individuals recognize that decisions made today affect those to be made in the future, at least in part, by expanding or contracting the set of admissible choices, that is, by lowering or raising the cost of a future choice. Such intertemporal linkages reside at the core of all dynamic processes in economics. Consequently, mathematical methods that account for such intertemporal linkages are fundamental, 
in principle, to all economic decisions''.  In the same book an economic 
interpretation of control problems and of Hamilton-Jacobi-Bellman equations are 
provided.

In \cite{hv}, the problem $(P)$ has been studied in the case without discount:  
in order to have the Hopf-Lax formula the Lagrangian depends only on the 
derivative $\dot y$. In particular, it does not depend on time.  The idea 
in this paper and in \cite{hv} is to embed the vector-valued problem into a 
set-valued one.  More precisely, it is fundamental to work in a complete lattice, 
so that infima and suprema are well defined.  For this approach see 
\cite{SetOptSurvey}.

In \cite{RZ}, the author proposed problem $(P)$ for a real-valued utility 
function, with a non-constant discount, provided a Hop-Lax formula and deduced 
a dynamic programming equation.

The present paper generalizes both \cite{hv} and \cite{RZ}. With respect to the 
first one, there is the non-constant discount. With respect to the second 
one, there is the multi-objective Lagrangian. To the knowledge of the author, 
this kind of generalization has never been addressed.

When the problem $(P)$ is embedded into a set-valued framework, a complete 
lattice structure is obtained and so the value function has a straightforward 
definition.  It is important to notice that, because of the discount factor, it can 
be called a current value optimal value function, for its value is discounted back 
to time $t$.  

The Hopf-Lax formula was found in the 50's (see \cite{Lax57} and \cite{Hopf65}, 
the first one in dimension one and the second one in general dimension).  This 
result came from the fact that straight lines are optimal trajectories.  This 
happens when the Hamiltonian function only depends on the derivative of the 
arc ($H=H(p)$).  In \cite{BJL} it is proved that this happens also when the 
Hamiltonian depends also on the arc $H=H(y,p)$, it is nondecreasing in $y$ and 
convex and positively homogeneous of degree 1 in $p$.

More recent developments can be found in \cite{Stroemberg}, \cite{Aubin}, \cite{AvantaggiatiLoreti}, \cite{Hoang}.

In \cite{RZ} as well as in this paper, the optimal trajectories are not linear, 
but all the same it is possible to write a generalization of the Hopf-Lax formula, in 
the sense that it is still possible to shift the infimum from an infinite dimensional 
space to a finite dimensional one.

As regards the importance of considering a non-constant discount factor, this 
choice arises in many economic models.  The utility function $L$, that measures 
the satisfaction of an agent, changes during the time, in the sense that an 
earlier attainment of the utility gives the agent a higher satisfaction.  This is 
why the discount factor is also called impatience rate (see for example 
\cite{Barro}, \cite{Karp}, \cite{MSolano_Navas}).

The problem of \cite{RZ} has been generalized to the field of deterministic
differential games with a stochastic terminal time in \cite{M-SS}.  See also 
\cite{M_SP}.

\section{Preliminaries}

First, some basic concepts and definitions of set optimization and of 
complete lattice approach are recalled.  For more on the subject, see \cite{SetOptSurvey}.

Let $C$ be a closed and convex cone in $\RR^d$ with nonempty interior. Its dual 
cone is defined as
$$
C^+=\{\zeta\in\RR^d\mid\forall z\in C\ \zeta\cdot z\geq 0\},
$$
where $\cdot$ denotes the scalar product in $\RR^d$.

If $\P(\RR^d)$ is the power set of $\RR^d$, it can be endowed with the 
Minkowski sum, i.e., for $A,B\in\P(\RR^d)$ one sets $A+B=\{ a+b\mid a\in A,\ 
b\in B\}$.  For the empty set, the sum is defined as $A+\emptyset=\emptyset+A
=\emptyset$.  
We shall consider also the closure with respect to the ususal topology of the sum:
$$
A\oplus B=\mbox{cl}\;(A+B).
$$

We denote by $\F(\RR^d,C)$ the subset of $\P(\RR^d)$ of those sets which are 
invariant with respect to the sum of the cone and closed:
$$
\F(\RR^d,C)=\{ A\subseteq\RR^d\mid A=A\oplus C\}.
$$
The pair $(\F(\RR^d,C),\supseteq)$ is a complete lattice, where the infimum 
and the supremum over a collection $\A\subseteq\F(\RR^d,C)$ can be found as
$$
\inf\A=\cl\bigcup_{A\in\A} A,\qquad\sup\A=\bigcap_{A\in\A} A.
$$

For any $\zeta\in C^+\bs\{0\}$, we define the half-space
$$
H^+(\zeta)=\{ z\in\RR^d\mid\zeta\cdot z\geq 0\}.
$$
We can also consider the $\zeta$-difference of two sets $A,B\in\F(\RR^d,C)$:
$$
A-_\zeta B=\{ z\in\RR^d\mid z+B\subseteq A\oplus H^+(\zeta)\}.
$$
It is possible to see that the difference can also be written as
$$
A-_\zeta B=\{ z\in\RR^d\mid \zeta\cdot z+\inf_{b\in B}\zeta\cdot b\geq
  \inf_{a\in A}\zeta\cdot a\}
$$
and that it is always a closed half-space or the empty set or the whole space.

We give here a definition of limit.  Let $\{ A_m\}_{m\in\NN}$ be a sequence of 
sets in $\F(\RR^d,C)$.  The notation $\lim_{m\to\infty}A_m$ identifies the set
$$
\lim_{m\to\infty}A_m=\{ z\in\RR^d\mid\forall m\in\NN,\ \exists z_m\in A_m
  \mbox{ and }\lim_{m\to\infty}z_m=z\},
$$
which is still in $\F(\RR^d,C)$.  This definition coincides with 
Painlev\'e-Kuratowski upper limit (see \cite{Aubin_Frank}).   Generalizing, if 
$\{ A_s\}_{s\in S}$ with $S\subseteq\RR$ is a family of sets in $\F(\RR^d,C)$ and 
$s_*\in\RR$, we denote by $\lim_{s\to s_*}A_s$ the set which satisfies for any 
sequence $\{ s_m\}_{m\in\NN}\subseteq S$, with $\lim_{m\to\infty}s_m=s_*$,
$$
\lim_{s\to s_*}A_s=\lim_{m\to\infty} A_{s_m}.
$$

Let $X$ be a vector space and $f \colon X\to \F(\RR^d, C)$ be a 
set-valued function. The graph of $f$ is the set
\[
\mbox{graph}\, f=\{(x,z)\in X\times\RR^d\mid z\in f(x)\} \subseteq X\times \RR^d
\]
and its (effective) domain is the set
\[
\mbox{dom}\, f=\{ x\in X\mid f(x)\neq\emptyset\}\subseteq X\, .
\]

Solutions of set optimization problems split into sets generating the infimum 
(infimizers) and points with a minimal function value with respect to the order in 
the corresponding lattice.  The solution concept is due to \cite{HL} (see also 
\cite{SetOptSurvey} for more comments and references).
For $M \subseteq X$, we set
\[
f[M]=\{f(x) \mid x \in M\} \, .
\]
The set $M \subseteq \mathcal X$ is called an infimizer of $f$ if
\begin{equation}
\label{Infimizer}
\inf f[M] = \inf f[X].
\end{equation}

We want now to define a set-valued function, which generalizes a linear function.  
Let $\eta\in\RR^n$ and $\zeta \in C^+$ be given. We consider 
$S_{(\eta,\zeta)}:\RR^n\to\F(\RR^d,C)$:
\begin{equation}
\label{S}
S_{(\eta,\zeta)}(x)=\{z\in\RR^d\ \vert\ \zeta\cdot z\geq\eta\cdot x\}\, .
\end{equation}
In fact this function is half-space valued and is additive and positively homogeneous 
in $x$.

The Fenchel conjugate of the function $f:\RR^n\to\F(\RR^d,C)$ is 
defined as the function
\begin{equation}
\label{Fenchel}
\begin{array}{cccc}
f^*: & \RR^n\times C^+\backslash\{0\} & \to & \F(\RR^d,C) \\
       & (\eta,\zeta) & \mapsto & \sup_{x\in\RR^n} S_{(\eta,\zeta)}(x)-_{z^*}f(x)
\end{array}
\end{equation}
%

\section{Variational problem with discount}

Let $L:\RR^n\to\RR^d$, $g:\RR^n\to\RR^d$ be two functions mapping into 
$\RR^d$.   Fix $T>0$ and consider a variable discount factor $d_t(s)$ for $t$ in 
$[0,T]$ and $s$ in $[t,T]$, taking positive values.
 
Let $\mathcal{L}_{t}:[0,T]\times\RR^n\to\F(\RR^d,C)$ denote the 
set-valued function
$$
\mathcal{L}_{t}(s,w)=d_t(s)L(w)+C\, .
$$

We define the functionals $I_t: W^{1,1}([t,T],\RR^n)\to\RR^d$, 
$J_t: W^{1,1}([t,T],\RR^n)\to\F(\RR^d,C)$ by
$$
\begin{aligned}
I_t(y) &= \int_t^Td_{t}(s)L(\dot y(s))\, ds+d_t(T)g(y(T))\, , \\
J_t(y) &= \int_t^T\mathcal{L}_{t}(s,\dot y(s))\, ds+d_t(T)g(y(T))\, ,
\end{aligned}
$$
where the second integral is understood in the Aumann sense (see \cite{Aumann}) 
and where for every $F:[0,T]\to\F(\RR^d,C)$, $t\in[0,T]$, we define
$$
\int_t^tF(s)\, ds=C\, .
$$

For any $x\in\RR^n$ we shall consider the problem:
\begin{equation}
\label{pb}
\begin{aligned}
\mbox{minimize } &J_t(y) \\
\mbox{over the set }& A(t,x)=\{ y\in W^{1,1}([t,T],\RR^n)\ \vert\ y(t)=x\}\, .
\end{aligned}
\end{equation}

Since the problem has been now embedded into a set-valued problem, we 
are now working on the complete lattice $\F(\RR^d,C)$.  This means that the 
infimum and the supremum are well defined and the value function is 
(see also \cite{hv}):
\begin{equation}
\label{valuef}
U(t,x)=\inf_{y\in A(t,x)} J_t(y)\, .
\end{equation}

For simplicity for any $\zeta\in C^+\bs\{0\}$, we denote
$$
\begin{aligned}
L_\zeta : &\RR^n \longrightarrow\RR \\
&w \longmapsto L(w)\cdot\zeta \\
g_\zeta : &\RR^n \longrightarrow\RR \\
&w \longmapsto g(w)\cdot\zeta \\
I_{t,\zeta} : &W^{1,1}([t,T],\RR^n) \longrightarrow\RR \\
&y \longmapsto I_t(y)\cdot\zeta
\end{aligned}
$$

Let $B^+$ be a base of $C^+$, i.e., for each element $\zeta\in C^+\bs\{0\}$ 
there exist a unique $\xi\in B^+$ and a unique $\lambda>0$ such that $\zeta=
\lambda\xi$.  For example $B^+$ can be formed by the unitary vectors in $C^+$.  
Another possibility is, if there is an element $c_0\in C$ such that $\zeta\cdot c_0>0$ 
for all $\zeta\in C^+\bs\{0\}$, then the set $B^+_{c_0}=\{\zeta\in C^+\mid\zeta
\cdot c_0=1\}$ is a base of $C^+$.

We will consider the following hypotheses:
\begin{itemize}
\item[(h1)] All the scalarizations $L_\zeta$, $\zeta\in B^+$ of $L$ are $C^2$, strictly 
  convex and coercive
  $$
  \lim_{|w|_n\to\infty}\frac{L_\zeta(w)}{|w|_n}=+\infty\, ,
  $$
  where $|\cdot|_n$ is the standard norm in $\RR^n$.
\item[(h2)] All the scalarizations $g_\zeta$, $\zeta\in B^+$ of $g$ are globally 
  Lipschitz in $\RR^n$.
\item[(h3)] For any $t\in[0,T]$, $d_t:[t,T]\to(a,1]$, with $a>0$, is continuous 
  and $d_t(t)=1$ for each $t\in[0,T]$.
\end{itemize}

Throughout the paper, $\nabla$ denotes the gradient of a real function or the 
Jacobian matrix of a vector function and $\nabla^2$ the Hessian matrix of a 
real function.

\begin{remark}
\label{rmk-lmm3.1}
  By Lemma 3.1 in \cite{RZ}  (see also \cite{DGZ}) and hypotheses (h1) and (h3), 
  for any $\zeta\in C^+\setminus\{0\}$ and fixed $s\in(t,T]$,
  \begin{itemize}
  \item[(i)]  $\nabla L_\zeta$ is a homeomorphism of $\RR^n$;
  \item[(ii)] the mappings from $\RR^n$ to $\RR^n$
    $$
    \begin{aligned} 
    p &\longmapsto (\nabla L_\zeta)^{-1}\left( \frac{p}{d_t(s)} \right)\\
    p &\longmapsto \int_t^s  (\nabla L_\zeta)^{-1}\left( \frac{p}{d_t(r)} \right)dr
    \end{aligned}
    $$
    are of class $C^1$ and surjective.
  \end{itemize}
\end{remark}

For any $p\in\RR^n$,  $\zeta\in C^+\bs\{0\}$, we consider the arc
\begin{equation}
\label{Ypzeta}
Y_{t,x,p,\zeta}(s)=x+\int_t^s  (\nabla L_\zeta)^{-1}\left( \frac{p}{d_t(r)} \right)dr\, .
\end{equation}
It is an element of $A(t,x)$.

In the following lemma we study some property of concavity of the Hamiltonian 
function, recalling the definition of the function $S$ in \eqref{S}.

\begin{lemma}
\label{lmm_Hconcave}
  Let $\zeta\in C_+\setminus\{0\}$, $(t,x)\in[0,T]\times\RR^n$ and 
  $\mathcal{H}_{t,\zeta}:[t,T]\times\RR^n\times\RR^n\to\F(\RR^d,C_+)$ 
  be defined by
  $$
  \mathcal{H}_{t,\zeta}(s,w,p)=S_{(p,\zeta)}(w)-_\zeta\mathcal{L}_t(s,w)\, .
  $$
  Then
  $$
  \mathcal{L}_t^*(s,p,\zeta)=\sup_{w\in\RR^n}\mathcal{H}_{t,\zeta}(s,w,p)=
    \mathcal{H}_{t,\zeta}\left(s,\dot Y_{t,x,p,\zeta}(s),p\right)\, .
  $$
\end{lemma}

\begin{proof}
  It is immediate to see that
  $$
  \mathcal{L}_t^*(s,p,\zeta)\subseteq\mathcal{H}_{t,\zeta}\left(s,\dot Y_{t,x,p,
    \zeta}(s),p\right)\, .
  $$
  For the other inclusion we notice that
  $$
  z\in\mathcal{H}_{t,\zeta}\left(s,w,p\right)
  $$
  if and only if
  $$
  z\cdot\zeta+d_t(s)L_\zeta(w)\geq p\cdot w\, .
  $$
  Since the real-valued function $p\cdot w-d_t(s)L_\zeta(w)$ is concave in $w$ and 
  has a maximizer in $\dot Y_{t,x,p,\zeta}(s)$, we have
  $$
  p\cdot \dot Y_{t,x,p,\zeta}(s)-d_t(s)L_\zeta\left(\dot Y_{t,x,p,\zeta}(s)
    \right)\geq p\cdot w-d_t(s)L_\zeta(w)
  $$
  for any $w\in\RR^n$.  Now, if $z\in\mathcal{H}_{t,\zeta}\left(s,\dot Y_{t,x,
  p,\zeta}(s),p\right)$,
  $$
  z\cdot\zeta\geq p\cdot \dot Y_{t,x,p,\zeta}(s)-d_t(s)L_\zeta\left(\dot 
  Y_{t,x,p,\zeta}(s)\right)\geq p\cdot w-d_t(s)L_\zeta(w)
  $$
  for any $w\in\RR^n$ and
  $$
  z\in\bigcap_{w\in\RR^n}\mathcal{H}_{t,\zeta}(s,w,p)\, .
  $$
\end{proof}

As a consequence of the previous lemma, the following property holds.

\begin{lemma}
\label{lmm_LinC}
  For any $(t,x)\in[0,T]\times\RR^n$, $w,p\in\RR^n$,  $\zeta\in C^+\bs
  \{0\}$ and $s\in[t,T]$, there holds
  $$
  p\cdot\left(\dot Y_{t,x,p,\zeta}(s)-w\right)\frac{\zeta}{|\zeta|_d^2}-
    d_t(s)\left(L\left(\dot Y_{t,x,p,\zeta}(s)\right)-L(w)\right)\in H^+(\zeta)\, ,
  $$
  where $|\cdot|_d$ is the norm in $\RR^d$.
\end{lemma}

\begin{proof}
  By Lemma~\ref{lmm_Hconcave}
  \begin{equation}
  \label{incl}
  S_{(p,\zeta)}\left(\dot Y_{t,x,p,\zeta}(s)\right)-_\zeta\mathcal{L}_t\left(s,\dot 
    Y_{t,x,p,\zeta}(s)\right)\subseteq S_{(p,\zeta)}(w)-_\zeta\mathcal{L}_t(s,w)\, .
  \end{equation}
  These sets are half-spaces because they are $\zeta$-differences and can be 
  written in the following way:
  $$
  S_{(p,\zeta)}(w)-_\zeta\mathcal{L}_t(s,w)=p\cdot w\frac{\zeta}{|\zeta|_d^2}-
    d_t(s)L(w)+H^+(\zeta)\, .
  $$
  Then the inclusion (\ref{incl}) can be written
  $$
  p\cdot\left(\dot Y_{t,x,p,\zeta}(s)- w\right)\frac{\zeta}{|\zeta|_d^2}-
    d_t(s)\left(L\left(\dot Y_{t,x,p,\zeta}(s)\right)-L(w)\right)\in H^+(\zeta)\, .
  $$
\end{proof}

\section{Bellman's optimality principle}

Usually, in Bellman's optimality principle an inequality and an equation are involved 
that link the value  function evaluated at two different times $t<\tau$.  
This is true for example when the Lagrangian is real-valued (see Lemma 
4.1 and Corollary 4.1 in \cite{RZ}).  Instead, if the Lagrangian is vector-valued, 
the situation is more complex and 
it is not possible to obtain the value function at time $\tau$.  In fact the infimum 
is taken over the sum of two parts, as one can see in the following theorem.

\begin{theorem}
  For every initial condition $(t,x)\in[0,T]\times\RR^n$, admissible arc $y\in A(t,x)$ 
  and $\tau\in [t,T]$, we have
  \begin{equation}
  \label{Bell}
  U(t,x)\supseteq \int_t^\tau \mathcal{L}_t(s,\dot y(s))\, ds\oplus
    \inf_{\eta\in A(\tau,y(\tau))} \left[ W(t,\tau,\eta)+J_\tau(\eta) \right],
  \end{equation}
  where $W:[0,T]\times[t,T]\times W^{1,1}([\tau,T],\RR^n)\to\RR^d$,
  $$
  W(t,\tau,\eta)=\int_\tau^T(d_t(s)-d_\tau(s))
    L(\dot\eta(s))\, ds+(d_t(T)-d_\tau(T))g(\eta(T))\, .
  $$
  Moreover, the set $M\subset A(t,x)$ is an infimizer for problem (\ref{pb}) (see 
  definition \eqref{Infimizer}) if and only if for all $\tau\in[t,T]$
  \begin{equation}
  \label{Bell1}
  U(t,x) = \inf_{y\in M} \left[ \int_t^\tau \mathcal{L}_t(s,\dot y(s))\, ds\oplus
    \inf_{\eta\in A(\tau,y(\tau))} \left[ W(t,\tau,\eta)+J_\tau(\eta) \right]\right] .
  \end{equation}
\end{theorem}

\begin{proof}
  Consider $y\in A(t,x)$.  If $\tau=T$, $A(T,y(T))=\{ y(T)\}$ and
  \begin{equation}
  \label{tau=T}
  \inf_{\eta\in A(T,y(T))} \left[ W(t,T,\eta)+J_T(\eta) \right]=
    (d_t(T)-1)g(y(T))+g(y(T))+C=d_t(T)g(y(T))+C
  \end{equation}
  and
  $$
  \int_t^T \mathcal{L}_t(s,\dot y(s))\, ds\oplus d_t(T)g(y(T))+C\subseteq U(t,x)\, .
  $$
  Also for $\tau =t$ the inclusion is trivially true.  Let us now consider $t<\tau<T$.  
  For any $\eta\in A(\tau,y(\tau))$, we can consider the arc
  $$
  y_\eta(s)=\left\{\begin{array}{ll}
    y(s) & \mbox{if }t\leq s\leq\tau\\
    \eta(s) & \mbox{if }\tau<s\leq T\end{array}\right.
  $$
  and $y_\eta\in A(t,x)$.  It is possible to write $J_t(y_\eta)$ as
  \begin{equation}
  \label{Jt(yeta)}
  \begin{aligned}
  J_t(y_\eta) &= \int_t^T\mathcal{L}_t(s,\dot y_\eta(s))ds+d_t(T)g(y_\eta(T)) \\
    &=\int_t^\tau\mathcal{L}_t(s,\dot y(s))ds\oplus\int_\tau^T\mathcal{L}_t
      (s,\dot\eta(s))ds+d_t(T)g(\eta(T)) \\
    &=\int_t^\tau\mathcal{L}_t(s,\dot y(s))ds\oplus [W(t,\tau,\eta)+J_\tau(\eta)]\, .
  \end{aligned}
  \end{equation}
  Since
  $$
  U(t,x)\supseteq J_t(y_\eta)
  $$
  for every $\eta$, then, using \eqref{Jt(yeta)}, one gets
  $$
  \begin{aligned}
  U(t,x) &\supseteq\inf_{\eta\in A(\tau,y(\tau))} \left(\int_t^\tau\mathcal{L}_t(s,
      \dot y(s))ds\oplus [W(t,\tau,\eta)+J_\tau(\eta)]\right) \\
    &=\int_t^\tau\mathcal{L}_t(s,\dot y(s))ds\oplus \inf_{\eta\in A(\tau,y(\tau))}
      [W(t,\tau,\eta)+J_\tau(\eta)]
  \end{aligned}
  $$
  and \eqref{Bell} is proved.
  
  If for all $\tau\in[t,T]$ (\ref{Bell1}) holds, in particular for $\tau=T$, using 
  (\ref{tau=T}),
  $$
  U(t,x)= \inf_{y\in M} \left[ \int_t^T \mathcal{L}_t(s,\dot y(s))\, ds\oplus 
    d_t(T)g(y(T))+C\right]= \inf_{y\in M} J_t(y)
  $$
  and $M$ is an infimizer.
  
  Finally, we want to prove that, if $M$ is an infimizer, then for any $\tau\in[t,T]$
  $$
  U(t,x)\subseteq \inf_{y\in M} \left[ \int_t^\tau \mathcal{L}_t(s,\dot y(s))\, ds
    \oplus\inf_{\eta\in A(\tau,y(\tau))} \left[ W(t,\tau,\eta)+J_\tau(\eta) 
    \right]\right] .
  $$
  It is sufficient to prove that for any $y\in M$
  $$
  J_t(y)\subseteq\inf_{y\in M} \left[ \int_t^\tau \mathcal{L}_t(s,\dot y(s))\, ds\oplus
    \inf_{\eta\in A(\tau,y(\tau))} \left[ W(t,\tau,\eta)+J_\tau(\eta) \right]\right] .
  $$
  We can write $J_t(y)$ in a similar way to (\ref{Jt(yeta)}), using $\eta(s)=y(s)$, 
  and this concludes the proof.
\end{proof}

\section{Hopf-Lax formula}

Before stating the Hopf-Lax formula, we need the following lemmas.   The 
first one is a coercivity result.

\begin{lemma}
\label{lmm_coerc}
  Given $\zeta\in C^+\setminus\{0\}$ and $(t,x)\in[0,T]\times\RR^n$, the 
  following limits hold:
  \begin{itemize}
  \item[(i)] $\lim_{|w|_n\to +\infty}|\nabla L_\zeta(w)|_n=+\infty$,
  \item[(ii)] $\lim_{|p|_n\to +\infty}|(\nabla L_\zeta)^{-1}(p)|_n=+\infty$,
  \item[(iii)] for any $s\in[t,T]$, $\lim_{|p|_n\to +\infty} \left| \dot Y_{t,x,p,\zeta}(s) 
    \right|_n =+\infty$,
  \item[(iv)] $\lim_{|p|_n\to +\infty} I_{t,\zeta}(Y_{t,x,p,\zeta})=+\infty$.
  \end{itemize}
\end{lemma}

\begin{proof}
  We want to prove that the functions $|\nabla L_\zeta(w)|_n$ and 
  $|(\nabla L_\zeta)^{-1}(p)|_n$ are coercive.  Since $L_\zeta$ is convex, there 
  holds
  $$
  L_\zeta(0)-L_\zeta(w)\geq\nabla L_\zeta(w)\cdot(-w)
  $$
  for any $w\in\RR^n$ and this gives that
  \begin{equation}
  \label{nabla1}
  |\nabla L_\zeta(w)|_n|w|_n\geq\nabla L_\zeta(w)\cdot w\geq L_\zeta(w)-L_\zeta(0)\, .
  \end{equation}
  Hypothesis (h1) implies that for $w\in\RR^n$, with $|w|_n$ sufficiently big, there 
  exist $K_\zeta,\epsilon_\zeta>0$ such that
  \begin{equation}
  \label{nabla2}
  L_\zeta(w)\geq K_\zeta|w|_n^{1+\epsilon_\zeta}\, .
  \end{equation}
  From (\ref{nabla1}) and (\ref{nabla2}), dividing by $|w|_n$, the following 
  inequality is obtained
  $$
  |\nabla L_\zeta(w)|_n\geq K_\zeta|w|_n^{\epsilon_\zeta}-\frac{L_\zeta(0)}{|w|_n}
  $$
  and, taking the limit, one obtains (i).
  
  For the inverse function let us assume that there exists a sequence $\{ p_m
  \}_{m\in\NN}\subseteq\RR^n$ such that
  \begin{equation}
  \label{pinfty}
  \lim_{m\to +\infty} |p_m|_n=+\infty
  \end{equation}
  and
  $$
  \lim_{m\to +\infty}|(\nabla L_\zeta)^{-1}(p_m)|_n=\alpha\, .
  $$
  Then the sequence $\{ w_m\}_{m\in\NN}\subseteq\RR^n$, defined by
  $$
  w_m=(\nabla L_\zeta)^{-1}(p_m)
  $$
  is bounded and there exists a converging subsequence (that we still denote 
  $\{ w_m\}_{m\in\NN}$) $w_m\to\hat w$.  Then we obtain that
  $$
  \lim_{m\to +\infty}|p_m|_n=\lim_{m\to +\infty} |\nabla L_\zeta(w_m)|_n=
    |\nabla L_\zeta(\hat w)|_n\, .
  $$
  This contradicts (\ref{pinfty}) and so (ii) holds.
  This implies in particular (iii).
  
  By hypothesis (h2) there exists a Lipschitz constant $G_\zeta>0$ for the 
  function $g_\zeta$:
  $$
  g_\zeta(0)-g_\zeta(Y_{t,x,p,\zeta}(T))\leq |g_\zeta(Y_{t,x,p,\zeta}(T))-g_\zeta(0)|
    \leq G_\zeta|Y_{t,x,p,\zeta}(T)|_n\, .  
  $$
  In order to find the coercivity with respect to $p$ of the following function, we 
  calculate for $|p|_n$ sufficiently big
  $$
  \begin{aligned}
  I_{t,\zeta}(Y_{t,x,p,\zeta}) = &\int_t^T d_t(s)L_\zeta\left(\dot Y_{t,x,p,\zeta}
      (s)\right)ds+d_t(T)g_\zeta(Y_{t,x,p,\zeta}(T)) \\
    \geq &\int_t^T d_t(s)K_\zeta\left|\dot Y_{t,x,p,\zeta}(s)\right|_n^{1+
      \epsilon_\zeta}ds+d_t(T)g_\zeta(0)-d_t(T)G_\zeta|Y_{t,x,p,\zeta}(T)|_n \\
    \geq &\int_t^T d_t(s)K_\zeta\left|(\nabla L_\zeta)^{-1}\left( \frac{p}{d_t(s)} \right)
      \right|_n^{1+\epsilon_\zeta}ds+d_t(T)g_\zeta(0) \\
    &-d_t(T)G_\zeta \left( |x|_n+\left| \int_t^T  (\nabla L_\zeta)^{-1}\left( 
      \frac{p}{d_t(s)}\right)ds \right|_n\right) .
  \end{aligned}
  $$
  Since, if $F\in L^1([a,b];\RR^n)$, then
  $$
  \left|\int_a^bF(s)\, ds\right|_n\leq\sqrt{n}\int_a^b|F(s)|_n\, ds\, ,
  $$
  we have that
  $$
  \begin{aligned}
  I_{t,\zeta}(Y_{t,x,p,\zeta}) \geq &\int_t^T d_t(s)K_\zeta\left|(\nabla L_\zeta)^{-1}
      \left( \frac{p}{d_t(s)} \right)\right|_n^{1+\epsilon_\zeta}ds+d_t(T)g_\zeta(0) \\
    &-d_t(T)G_\zeta \left( |x|_n+\sqrt{n}\int_t^T \left|  (\nabla L_\zeta)^{-1}
      \left( \frac{p}{d_t(s)}\right)\right|_n ds \right) \\
    = & \int_t^T \left( d_t(s)K_\zeta\left|(\nabla L_\zeta)^{-1}
      \left( \frac{p}{d_t(s)} \right)\right|_n^{\epsilon_\zeta}-d_t(T)\sqrt{n}G_\zeta 
      \right) \left|(\nabla L_\zeta)^{-1}\left( \frac{p}{d_t(s)} \right)\right|_nds \\
    &+d_t(T)g_\zeta(0) -d_t(T)G_\zeta |x|_n
  \end{aligned}
  $$
  and this proves (iv).
\end{proof}

In the following lemma, an arc is given, that minimizes the functional with respect to 
every direction of the dual cone.

\begin{lemma}
\label{lmm-minp}
  Given $\zeta\in C^+\setminus\{0\}$ and $(t,x)\in[0,T]\times\RR^n$, there 
  exists $p(t,x,\zeta)\in\RR^n$ such that
  \begin{equation}
  \label{p(t,x,zeta)}
  \begin{aligned}
  \inf_{p\in\RR^n}I_{t,\zeta}(Y_{t,x,p,\zeta})
    &=I_{t,\zeta}(Y_{t,x,p(t,x,\zeta),\zeta})\, , \\
  \inf_{p\in\RR^n} J_t\left(Y_{t,x,p,\zeta}\right)+ H^+(\zeta)
    &= J_t(Y_{t,x,p(t,x,\zeta),\zeta})+ H^+(\zeta)\, .
  \end{aligned}
  \end{equation}
\end{lemma}

We observe that
$$
\inf_{p\in\RR^n} \left[ J_t\left(Y_{t,x,p,\zeta}\right) \right]+ H^+(\zeta) = 
  \inf_{p\in\RR^n} \left[ J_t\left(Y_{t,x,p,\zeta}\right)+ H^+(\zeta) \right].
$$

\begin{proof}
  By (iv) of Lemma \ref{lmm_coerc} the function of $p$ $I_{t,\zeta}
  (Y_{t,x,p,\zeta})$ attains its minimum at $p(t,x,\zeta)=p_0\in\RR^n$.

  It is obvious that
  $$
  \inf_{p\in\RR^n} J_t\left(Y_{t,x,p,\zeta}\right)+ H^+(\zeta)\supseteq 
    J_t(Y_{t,x,p_0,\zeta})+ H^+(\zeta)\, .
  $$
  In order to prove that
  $$
  \inf_{p\in\RR^n} J_t\left(Y_{t,x,p,\zeta}\right)+ H^+(\zeta)\subseteq 
    J_t(Y_{t,x,p_0,\zeta})+ H^+(\zeta)\, ,
  $$
  we consider $z\in J_t\left(Y_{t,x,p,\zeta}\right)+H^+(\zeta)$ for some 
  $p\in\RR^n$, then we can write
  $$
  z= I_{t,\zeta}\left(Y_{t,x,p,\zeta}\right) \frac{\zeta}{|\zeta|_d^2}+h
  $$
  with $h\in H^+(\zeta)$.  Since we have
  $$
  \begin{aligned}
  z &=  I_{t,\zeta}\left(Y_{t,x,p_0,\zeta}\right) \frac{\zeta}{|\zeta|_d^2} + 
      \left( I_{t,\zeta}\left(Y_{t,x,p,\zeta}\right) - I_{t,\zeta}\left(Y_{t,x,p_0,\zeta}
      \right) \right) \frac{\zeta}{|\zeta|_d^2} +h \\
    &= I_{t,\zeta}\left(Y_{t,x,p_0,\zeta}\right) \frac{\zeta}{|\zeta|_d^2} + h' 
      = I_t\left(Y_{t,x,p_0,\zeta}\right)+h''
  \end{aligned}
  $$
  with $h',h''\in H^+(\zeta)$, so $z\in J_t(Y_{t,x,p_0,\zeta})+ H^+(\zeta)$ 
  and this completes the proof.
\end{proof}

To simplify the notation we define now
\begin{equation}
Y_{t,x,\zeta}(s)=Y_{t,x,p(t,x,\zeta),\zeta}(s),
\end{equation}
where $p(t,x,\zeta)$ is defined in the previous lemma.

Now the Hopf-Lax formula can be stated.

\begin{theorem}
  Let $g$ have convex components.  If $x\in\RR^n$ and $0\leq t<T$, the value 
  function (\ref{valuef}) of problem (\ref{pb}) can be written as
  \begin{equation}
  \label{Hopf_Lax}
  U(t,x)=\sup_{\zeta\in C^+} \left( \inf_{p\in\RR^n} J_t(Y_{t,x,p,\zeta})+H^+(\zeta)
    \right) .
  \end{equation}
\end{theorem}

\begin{proof}
  Let us define
  \begin{equation}
  \label{V(t,x)}
  V(t,x) = \sup_{\zeta\in C^+} \left(\inf_{p\in\RR^n} J_t(Y_{t,x,p,\zeta})
    +H^+(\zeta)\right) = \sup_{\zeta\in C^+} \left( J_t(Y_{t,x,\zeta})+
    H^+(\zeta) \right),
  \end{equation}
  where Lemma \ref{lmm-minp} has been used. 
  We consider an arc $y\in A(t,x)$.  By Remark~\ref{rmk-lmm3.1} (ii), for any 
  $\zeta\in C^+\bs\{0\}$, there exists $\overline p$ such that $Y_{t,x,
  \overline p,\zeta}(T)=y(T)$.  By Lemma~\ref{lmm_LinC} there holds
  $$
  \overline p\cdot\left(\dot Y_{t,x,\overline p,\zeta}(s)-\dot y(s)\right)
    \frac{\zeta}{|\zeta|_d^2}-d_t(s)
    \left(L\left(\dot Y_{t,x,\overline p,\zeta}(s)\right)-L(\dot y(s))\right)\in H^+(\zeta)\, .
  $$
  This implies that
  $$
  \mathcal{L}_t(s,\dot y(s))\subseteq\mathcal{L}_t\left(s,\dot Y_{t,x,\overline p,\zeta}
    (s)\right)-\overline p\cdot\left(\dot Y_{t,x,\overline p,\zeta}(s)-\dot y(s)\right)
    \frac{\zeta}{|\zeta|_d^2}+H^+(\zeta)\, .
  $$
  Integrating the previous inclusion from $t$ to $T$, one obtains
  $$
  \begin{aligned}
  \int_t^T \mathcal{L}_t(s,\dot y(s))\, ds &\subseteq \int_t^T \mathcal{L}_t
      \left(s,\dot Y_{t,x,\overline p,\zeta}(s)\right)\, ds-\left(\int_t^T \overline p\cdot
      \left(\dot Y_{t,x,\overline p,\zeta}(s)-\dot y(s)\right)\, ds\right)
      \frac{\zeta}{|\zeta|_d^2}+H^+(\zeta) \\
    &=\int_t^T \mathcal{L}_t\left(s,\dot Y_{t,x,\overline p,\zeta}(s)\right)\, ds
      +H^+(\zeta)\, .
  \end{aligned}
  $$
  Adding $d_t(T)g(y(T))=d_t(T)g(Y_{t,x,\overline p,\zeta}(T))$ to both sets, we 
  obtain that for any arc $y$ there exists $\overline p$ such that
  $$
  \int_t^T \mathcal{L}_t(s,\dot y(s))\, ds+d_t(T)g(y(T))\subseteq
    \int_t^T \mathcal{L}_t\left(s,\dot Y_{t,x,\overline p,\zeta}(s)\right)\, ds+
    d_t(T)g(Y_{t,x,\overline p,\zeta}(T))+H^+(\zeta)
  $$
  and this proves that $U(t,x)\subseteq \inf_{p\in\RR^n} J_t(Y_{t,x,p,\zeta})
  +H^+(\zeta)$ for every $\zeta\in C^+\bs\{ 0\}$ and consequently
  $$
  U(t,x)\subseteq\sup_{\zeta\in C^+} \left(\inf_{p\in\RR^n} J_t(Y_{t,x,p,\zeta})
    +H^+(\zeta)\right) =V(t,x)\, .
  $$
  
  Let us suppose that there exists $z_0\in V(t,x)\bs U(t,x)$.  
  For every $\zeta\in C^+\bs\{0\}$,
  \begin{equation}
  \label{z0}
  z_0\in J_t(Y_{t,x,\zeta})+H^+(\zeta)
  \end{equation}
  and $z_0=I_t(Y_{t,x,\zeta})+h_\zeta$.  If we fix $c_0\in
  \mbox{int}\, C$ and consider the half-line $z_0+rc_0$ with $r>0$, for $r$ 
  sufficiently large $h_\zeta+rc_0\in C$ and consequently $z_0+rc_0\in U(t,x)$.  
  Then there exists $r_0>0$ such that $z_1=z_0+r_0c_0\in U(t,x)$ and it is on the 
  boundary of $U(t,x)$.  Since $U(t,x)$ is convex, by the supporting hyperplane 
  theorem (see for example \cite{BV}) there exists $\xi\in\RR^d\bs\{0\}$ 
  such that for any $z\in U(t,x)$
  $$
  \xi\cdot z\geq\xi\cdot z_1\, .
  $$
  From the fact that $z_1+c\in U(t,x)$ for any $c\in C$, we have that $\xi\cdot
  (z_1+c)\geq\xi\cdot z_1$ and $\xi\in C^+\bs\{0\}$.  From (\ref{z0}) 
  with $\zeta=\xi$, we obtain that
  \begin{equation}
  \label{xiz0}
  \xi\cdot z_0\geq I_{t,\xi}(Y_{t,x,\xi})\, .
  \end{equation}
  Since $I_t(Y_{t,x,\xi})\in U(t,x)$, we have
  $$
  I_{t,\xi}(Y_{t,x,\xi})\geq\xi\cdot z_1=\xi\cdot(z_0+r_0c_0)\, .
  $$
  This inequality and inequality (\ref{xiz0}) give
  $$
  r_0\xi\cdot c_0\leq 0\, ,
  $$
  but this implies that $\xi\cdot c_0=0$ and this is not possible because $c_0\in
  \mbox{int}\,C$.
  \end{proof}

\section{The Hamilton-Jacobi-Bellman equation}

In this section we assume that:
\begin{itemize}
\item[(h4)] the discount factor is of class $C^1$ in $t$ and 
\item[(h5)] all the scalarizations $g_\zeta$, with $\zeta$ in a base $B^+$ of 
  $C^+$ are $C^2$ and convex.
\end{itemize}

In the following lemma the differentiability of the arcs $Y_{t,x,p,\zeta}(s)$ and of 
their derivatives $\dot Y_{t,x,p,\zeta}(s)$ with respect to the parameters is studied.

\begin{lemma}
  Given $\zeta\in C^+\bs\{0\}$, $(t,x)\in[0,T]\times\RR^n$ and $p\in\RR^n$
  the arcs $Y_{t,x,p,\zeta}(s)$ and the derivatives $\dot Y_{t,x,p,\zeta}(s)$ admit 
  the partial derivatives with respect to $t$, $x$ and $p$:
  \begin{eqnarray*}
  \frac{\partial Y_{t,x,p,\zeta}(s)}{\partial t} &=& -\nabla L_\zeta^{-1}(p)
    -\int_t^s \frac{\partial d_t(r)}{\partial t} \nabla(\nabla L_\zeta)^{-1}\left(
    \frac{p}{d_t(r)}\right)\frac{p}{(d_t(r))^2}dr\, , \\
  \frac{\partial\dot Y_{t,x,p,\zeta}(s)}{\partial t} &=& -\frac{\partial d_t(s)}{\partial t}
    \nabla(\nabla L_\zeta)^{-1}\left(\frac{p}{d_t(s)}\right)\frac{p}{(d_t(s))^2}\, , \\
  \nabla_x Y_{t,x,p,\zeta}(s) &=& I\, , \\
  \nabla_x\dot Y_{t,x,p,\zeta}(s) &=& 0\, , \\
  \nabla_pY_{t,x,p,\zeta}(s) & = & \int_t^s \frac{1}{d_t(r)} \nabla
    (\nabla L_\zeta)^{-1} \left( \frac{p}{d_t(r)} \right) dr\, , \\
  \nabla_p\dot Y_{t,x,p,\zeta}(s) & = & \frac{1}{d_t(s)} \nabla
    (\nabla L_\zeta)^{-1} \left( \frac{p}{d_t(s)} \right)\, ,
  \end{eqnarray*}
  where $I$ is the identity in $\RR^n$ and $0$ is the null matrix.
\end{lemma}

\begin{proof}
  By Remark~\ref{rmk-lmm3.1} the arcs $Y_{t,x,p,\zeta}(s)$ are of class $C^1$ 
  in $p$.  We recall that $L_\zeta$ is of class $C^2$ and it is easy to see that
  $$
  \nabla(\nabla L_\zeta)^{-1}(x)=\left[\nabla^2L_\zeta((\nabla L_\zeta)^{-1}(x))
    \right]^{-1},
  $$
  so it exists and it is continuous.
\end{proof}

In the following proposition the function $p(t,x,\zeta)$ introduced in Lemma 
\ref{lmm-minp} is studied.

\begin{proposition}
  Given $\zeta\in C^+\bs\{0\}$ and $(t,x)\in[0,T]\times\RR^n$
  \begin{itemize}
  \item[(i)] $p(t,x,\zeta)$ in Lemma \ref{lmm-minp} is a solution of
    \begin{equation}
    \label{F=0}
    F(t,x,p,\zeta)=0
    \end{equation}
    where $F:[0,T]\times\RR^n\times\RR^n\times C^+\bs\{0\}\to\RR^n$ 
    is defined by
    $$
    F(t,x,p,\zeta)=p+d_t(T)\nabla g_\zeta\left(Y_{t,x,p,\zeta}(T)\right)
    $$
  \item[(ii)] The Jacobian matrix
    $$
    \nabla_pF (t,x,p,\zeta)=I+A(t,x,p,\zeta)\, ,
    $$
    where $I$ is the identity in $\RR^n$ and
    $$
    A(t,x,p,\zeta)=d_t(T)\nabla^2g_\zeta(Y_{t,x,p,\zeta}(T))
      \nabla_pY_{t,x,p,\zeta}(T)\, ,
    $$
    is nonsingular.
  \item[(iii)] The function $p(t,x,\zeta)$ is well defined in a neighborhood of $(t,x)$ 
    and admits the partial derivatives
    $$
    \begin{aligned}
    \frac{\partial p}{\partial t}(t,x,\zeta) = &-(I+A(t,x,p,\zeta))^{-1}\left[ 
        \frac{\partial d_t(T)}{\partial t}\nabla g_\zeta\left(Y_{t,x,p,\zeta}(T)
        \right) \right. \\
      &+ \left. d_t(T)\nabla^2g_\zeta(Y_{t,x,p,\zeta}(T)) 
        \frac{\partial Y_{t,x,p,\zeta}(T)}{\partial t}\right], \\
    \frac{\partial p}{\partial x}(t,x,\zeta) &= -(I+A(t,x,p,\zeta))^{-1}\left[ d_t(T)
      \nabla^2g_\zeta(Y_{t,x,p,\zeta}(T)) \nabla_x Y_{t,x,p,\zeta}(T)\right].
    \end{aligned}
    $$
  \end{itemize}
\end{proposition}

\begin{proof}
  From the definition in the first equation of \eqref{p(t,x,zeta)}, $p(t,x,\zeta)$ 
  must solve the equation
  $$
  \nabla_p I_{t,\zeta}(Y_{t,x,p,\zeta})=0\, .
  $$
  Calculating the previous derivative, one obtains
  $$
  \left[ p+d_t(T)\nabla g_\zeta(Y_{t,x,p,\zeta}(T)) \right]\cdot\int_t^T
    \frac{1}{d_t(s)}\nabla(\nabla L_\zeta)^{-1}\left(\frac{p}{d_t(s)}\right)ds=0\, .
  $$
  The solutions of equation \eqref{F=0} are obviously also solutions of the previous 
  equation.
  
  The matrix $A(t,x,p,\zeta)$ is the product of two matrices.  The first one is the 
  Hessian matrix of a $C^2$ function, so it is symmetric. The second one is the 
  integral of the matrix
  $$
  \frac{1}{d_t(r)}\nabla(\nabla L_\zeta)^{-1} \left( \frac{p}{d_t(r)} \right) = 
    \frac{1}{d_t(r)} \left[ \nabla^2L_\zeta \left( \nabla L_\zeta^{-1}\left( 
    \frac{p}{d_t(r)} \right)\right)\right]^{-1}.
  $$
  Since the Hessian matrix of $L_\zeta$ is symmetric, so it is its inverse and its 
  integral.  Moreover, the first matrix $\nabla^2 g_\zeta$ is positively semidefinite, 
  while the second one $\nabla_p Y_{t,x,p,\zeta}(T)$ is positively definite.  
  Then their product $A(t,x,p,\zeta)$ is also positively semidefinite (because they 
  can be simultaneously diagonalized).   
  If $\nabla_pF(t,x,p,\zeta)$ were singular, $A(t,x,p,\zeta)$ should have an 
  eigenvector of $-1$ and this is in contradiction with the fact that it is positively 
  semidefinite.  
\end{proof}

For $(t,x)\in[0,T]\times\RR^n$, $q\in\RR^n$ and $\zeta\in C^+\bs\{0\}$, we 
consider the partial derivatives:
\begin{equation}
\begin{aligned}
U_{t,\zeta}(t,x) &= \lim_{s\to 0^+} \frac{1}{s} \left[ U(t+s,x)-_\zeta U(t,x) \right],\\
U_{q,\zeta}(t,x) &= \lim_{s\to 0^+} \frac{1}{s} \left[ U(t,x+sq)-_\zeta
  U(t,x) \right].
\end{aligned}
\label{def_derivatives}
\end{equation}

Similar definitions are used in \cite{AC} and \cite{Pilecka}.   
These derivatives, if they exist, are closed half-spaces with normal $\zeta$ or in 
the extreme cases they are the empty set or $\RR^d$.

\begin{proposition}
  Given $\zeta\in C^+\bs\{0\}$ and $(t,x)\in[0,T]\times\RR^n$, the partial 
  derivatives with respect to $t$ and with respect to $x$ 
  in the direction $q$ exist and are the following ones:
  \begin{eqnarray}
  U_{t,\zeta}(t,x) &=& S_{(u_{t,\zeta}(t,x),\zeta)}(1) \label{Ut} \\
  U_{q,\zeta}(t,x) &=& S_{(\nabla u_\zeta(t,x),\zeta)}(q) \label{Ux}
  \end{eqnarray}
  where
  \begin{eqnarray*}
  u_{t,\zeta}(t,x) &=& -L_\zeta(\dot Y_{t,x,\zeta}(t))+\int_t^T
      \frac{\partial d_t(s)}{\partial t}L_\zeta(\dot Y_{t,x,\zeta}(s))ds 
      + \frac{\partial d_t(T)}{\partial t}g_\zeta(Y_{t,x,\zeta}(T)) \\
    &&- d_t(T)\nabla g_\zeta(Y_{t,x,\zeta}(T))\cdot
      (\nabla L_\zeta)^{-1}(p(t,x,\zeta))\, , \\
  \nabla u_\zeta(t,x) &=& d_t(T)\nabla g_\zeta(Y_{t,x,\zeta}(T))\, .
  \end{eqnarray*}
\end{proposition}

\begin{proof}
  First of all, we calculate, using Hopf-Lax formula and \eqref{p(t,x,zeta)},
  $$
  \begin{aligned}
  \inf\{\zeta\cdot z \mid z\in U(t,x)\} &=\inf \left\{\zeta\cdot z\mid z\in\sup_{\xi\in C^+}
      \left(\inf_{p\in\RR^n}J_t(Y_{t,x,p,\xi})+H^+(\xi)\right) \right\} \\
    &\geq \inf \left\{\zeta\cdot z\mid z\in\inf_{p\in\RR^n}J_t(Y_{t,x,p,\zeta})+
      H^+(\zeta) \right\} \\
    &= \inf\left\{ I_{t,\zeta}(Y_{t,x,p,\zeta})\mid p\in\RR^n\right\} =
      I_{t,\zeta}(Y_{t,x,\zeta}).
  \end{aligned}
  $$
  Since $I_t(Y_{t,x,\zeta})\in U(t,x)$, we can conclude that
  $$
  \inf\{\zeta\cdot z \mid z\in U(t,x)\} = I_{t,\zeta}(Y_{t,x,\zeta}).
  $$
  Now it is possible to write for $h>0$:
  $$
  \frac{1}{h}[U(t+h,x)-_\zeta U(t,x)]=\left\{ z\in\RR^d\mid \zeta\cdot z\geq\frac{1}{h}
    [I_{t+h,\zeta}(Y_{t+h,x,\zeta})-I_{t,\zeta}(Y_{t,x,\zeta})]\right\}.
  $$
  Since the total derivative with respect to the time $t$ is
  $$
  \begin{aligned}
  \lefteqn{\frac{d I_{t,\zeta}(Y_{t,x,\zeta})}{dt}=  -L_\zeta(\dot Y_{t,x,\zeta}(t)) 
      + \int_t^T \frac{\partial d_t(s)}{\partial t}L_\zeta(\dot Y_{t,x,\zeta}(s))ds 
      + \frac{\partial d_t(T)}{\partial t}g_\zeta(Y_{t,x,\zeta}(T))} \\
    &+ \int_t^T p(t,x,\zeta)\cdot\frac{\partial\dot Y_{t,x,\zeta}(s)}{\partial t}ds 
      + d_t(T)\nabla g_\zeta(Y_{t,x,\zeta}(T))\cdot
      \frac{\partial Y_{t,x,\zeta}(T)}{\partial t} \\
    &+ \int_t^T p(t,x,\zeta)\cdot \nabla_p\dot Y_{t,x,p(t,x,\zeta),\zeta}(s)ds
      \frac{\partial p(t,x,\zeta)}{\partial t} \\
    &+ d_t(T)\nabla g_\zeta(Y_{t,x,\zeta}(T))\cdot
      \nabla_p Y_{t,x,p(t,x,\zeta),\zeta}(T) \frac{\partial p(t,x,\zeta)}{\partial t} \\
    = &-L_\zeta(\dot Y_{t,x,\zeta}(t)) + \int_t^T
      \frac{\partial d_t(s)}{\partial t}L_\zeta(\dot Y_{t,x,\zeta}(s))ds \\
    &+ \frac{\partial d_t(T)}{\partial t}g_\zeta(Y_{t,x,\zeta}(T)) 
      -d_t(T)\nabla g_\zeta(Y_{t,x,\zeta}(T))\cdot\nabla L_\zeta^{-1}(p(t,x,\zeta)) \\
    &-\left(p(t,x,\zeta)+d_t(T)\nabla g_\zeta(Y_{t,x,\zeta}(T))\right)\cdot
      \int_t^T\frac{\partial d_t(s)}{\partial t}\nabla(\nabla L_\zeta)^{-1}
      \left(\frac{p(t,x,\zeta)}{d_t(s)}\right)\frac{p(t,x,\zeta)}{(d_t(s))^2}ds \\
    &+ \left(p(t,x,\zeta)+d_t(T)\nabla g_\zeta(Y_{t,x,\zeta}(T))\right)\cdot
      \int_t^T\nabla(\nabla L_\zeta)^{-1}\left(\frac{p(t,x,\zeta)}{d_t(s)}\right)
      \frac{1}{d_t(s)}ds\frac{\partial p(t,x,\zeta)}{\partial t} \\
    = &-L_\zeta(\dot Y_{t,x,\zeta}(t)) + \int_t^T
      \frac{\partial d_t(s)}{\partial t}L_\zeta(\dot Y_{t,x,\zeta}(s))ds \\
    &+ \frac{\partial d_t(T)}{\partial t}g_\zeta(Y_{t,x,\zeta}(T)) 
      -d_t(T)\nabla g_\zeta(Y_{t,x,\zeta}(T))\cdot\nabla L_\zeta^{-1}(p(t,x,\zeta))\, ,
  \end{aligned}
  $$
  the derivative with respect to time of the value function is as in \eqref{Ut}.
  
  On the other hand, we have
  $$
  \frac{1}{h}[U(t,x+hq)-_\zeta U(t,x)]=\left\{ z\in\RR^d\mid\zeta\cdot z\geq\frac1h
    [I_{t,\zeta}(Y_{t,x+hq,\zeta})-I_{t,\zeta}(Y_{t,x,\zeta})]
    \right\}.
  $$
  In order to prove \eqref{Ux}, we calculate the total derivative with respect to $x$
  $$
  \begin{aligned}
  \nabla_x &I_{t,\zeta}(Y_{t,x,\zeta}) = d_t(T)\nabla g_\zeta
      (Y_{t,x,\zeta}(T)) \\
    &+ (p(t,x,\zeta)+d_t(T)\nabla g_\zeta
      (Y_{t,x,\zeta}(T)))\cdot\int_t^T\nabla(\nabla L_\zeta)^{-1}\left(
      \frac{p(t,x,\zeta)}{d_t(s)}\right)\frac{1}{d_t(s)}ds D_xp(t,x,\zeta) \\
    =& d_t(T)\nabla g_\zeta(Y_{t,x,\zeta}(T))
  \end{aligned}
  $$
  and this concludes the proof.
\end{proof}

\begin{remark}
  It is not difficult to see that
  $$
  U(t+h,x+hq)-_\zeta U(t,x)=(U(t+h,x+hq)-_\zeta U(t+h,x))+(U(t+h,x)-_\zeta U(t,x))
  $$
  and that
  $$
  \lim_{h\to 0}\frac{1}{h}(U(t+h,x+hq)-_\zeta U(t+h,x))=U_{q,\zeta}(t,x).
  $$
  This means that
  $$
  \lim_{h\to0}\frac{1}{h}(U(t+h,x+hq)-_\zeta U(t,x))=U_{t,\zeta}(t,x)+U_{q,\zeta}(t,x).
  $$
\end{remark}

Recalling the Fenchel conjugate \eqref{Fenchel}, the following theorem can be 
stated.

\begin{theorem}
  Given $\zeta\in C^+\bs\{0\}$ and $(t,x)\in[0,T]\times\RR^n$, considering 
  $\mathcal{L}_t(t,\cdot):\RR^n\to\F(\RR^d,C)$, the value function $U(t,x)$ is 
  a solution of the equation
  \begin{equation}
  \label{Eq_HJB_zeta}
  U_{t,\zeta}(t,x)=\mathcal{L}_t^*(t,-\nabla u_\zeta(t,x),\zeta)+w(t,x,\zeta)\, ,
  \end{equation}
  where
  $$
  w(t,x,\zeta)=\int_t^T\frac{\partial d_t(s)}{\partial t} L\left(\dot Y_{t,x,\zeta}(s)
    \right)ds + \frac{\partial d_t(T)}{\partial t}g(Y_{t,x,\zeta}(T))\, .
  $$
\end{theorem}

\begin{proof}
  Given $q\in\RR^n$, we consider $y(s)=x+(s-t)q$, $y\in A(t,x)$.  Using Bellman's 
  inequality \eqref{Bell}, we have for $h>0$ sufficiently small
  $$
  \begin{aligned}
  U(t,x) &\supseteq\int_t^{t+h} \mathcal{L}_t(s,q)ds\oplus\inf_{\eta\in A(t+h,x+hq)}
      [W(t,t+h,\eta)+J_{t+h}(\eta)] \\
    &\supseteq\int_t^{t+h} \mathcal{L}_t(s,q)ds\oplus W(t,t+h,
      Y_{t+h,x+hq,\zeta})+J_{t+h}(Y_{t+h,x+hq,\zeta}).
  \end{aligned}
  $$
  Using now this inclusion of sets, we obtain
  $$
  \begin{aligned}
  U &(t+h,x+hq) -_\zeta U(t,x)\subseteq \\
    &U(t+h,x+hq)-_\zeta\left[ \int_t^{t+h} \mathcal{L}_t(s,q)ds\oplus W(t,t+h,
      Y_{t+h,x+hq,\zeta})+J_{t+h}(Y_{t+h,x+hq,\zeta}) \right].
  \end{aligned}
  $$
  We can calculate that
  $$
  \begin{aligned}
  \inf_{z\in U(t+h,x+hq)} \zeta\cdot z = &I_{t+h,\zeta}(Y_{t+h,x+hq,\zeta}) \\
  \inf_{z\in\int_t^{t+h} \mathcal{L}_t(s,q)ds} = &\int_t^{t+h} d_t(s)L_\zeta(q)ds \\
  \inf_{z\in W(t,t+h,Y_{t+h,x+hq,\zeta})+J_{t+h}
      (Y_{t+h,x+hq,\zeta})} \zeta\cdot z = &\zeta\cdot 
      W(t,t+h,Y_{t+h,x+hq,\zeta}) \\
    &+I_{t+h,\zeta}(Y_{t+h,x+hq,\zeta}).
  \end{aligned}
  $$
  Then we have the following inclusion:
  $$
  U(t+h,x+hq)-_\zeta U(t,x)\subseteq -\int_t^{t+h} d_t(s)L(q)ds-
    W(t,t+h,Y_{t+h,x+hq,\zeta})+H^+(\zeta)
  $$
  and, taking the limit, the partial derivatives fulfill the inclusion
  $$
  U_{t,\zeta}(t,x)+U_{q,\zeta}(t,x)\subseteq -L(q)+w(t,x,\zeta)+H^+(\zeta).
  $$
  As a consequence, we can write that
  $$
  U_{t,\zeta}(t,x)\subseteq (S_{(-\nabla u_\zeta(t,x),\zeta)}(q)-_\zeta L(q))+
    w(t,x,\zeta)
  $$
  for every $q\in\RR^n$.  It is then possible to take the supremum with respect to $q$
  \begin{equation}
  \label{Utsubset}
  \begin{aligned}
  U_{t,\zeta}(t,x) &\subseteq\sup_{q\in\RR^n}(S_{(-\nabla u_\zeta(t,x),\zeta)}(q)
      -_\zeta L(q))+w(t,x,\zeta) \\
    &=\mathcal{L}_t^*(t,-\nabla u_\zeta(t,x))+w(t,x,\zeta).
  \end{aligned}
  \end{equation}

  It is easy to check that for $q=\nabla L_\zeta^{-1}(p(t,x,\zeta))$
  $$
  U_{t,\zeta}(t,x)+U_{\nabla L_\zeta^{-1}(p),\zeta}(t,x)=-L\left(\nabla 
    L_\zeta^{-1}(p(t,x,\zeta))\right)+w(t,x,\zeta)+H^+(\zeta)
  $$
  and so
  $$
  U_{t,\zeta}(t,x)=\left[S_{(-\nabla u_\zeta(t,x),\zeta)}\left(\nabla L_\zeta^{-1}
    (p(t,x,\zeta))\right)-_\zeta L\left(\nabla L_\zeta^{-1}(p(t,x,\zeta))\right)\right]
    +w(t,x,\zeta).
  $$
  From \eqref{Utsubset} and the previous equation, one concludes that
  $$
  \begin{aligned}
  U_{t,\zeta}(t,x) &\subseteq\mathcal{L}_t^*(t,-\nabla u_\zeta(t,x))+
      w(t,x,\zeta) \\
    &\subseteq\left[S_{(-\nabla u_\zeta(t,x),\zeta)}\left(\nabla L_\zeta^{-1}
      (p(t,x,\zeta))\right)-_\zeta L\left(\nabla L_\zeta^{-1}(p(t,x,\zeta))\right)\right]
      +w(t,x,\zeta)=U_{t,\zeta}(t,x).
  \end{aligned}
  $$
  This proves \eqref{Eq_HJB_zeta}.
\end{proof}

In the following corollary, the Hamilton-Jacobi-Bellman equation is written 
independently of the directions in the dual cone.  It is easy to see that equation 
\eqref{Eq_HJB_zeta} can also be written
$$
U_{t,\zeta}(t,x)-_\zeta\mathcal{L}_t^*(t,-\nabla u_\zeta(t,x),\zeta)-w(t,x,\zeta)
  =H^+(\zeta)\, .
$$
However, the neutral element with respect to $\oplus$ in $\F(\RR^d,C)$ is $C$ and 
not $H^+(\zeta)$.  So, in order to have a ``complete'' equation and not one that 
describes only one direction, an intersection of the corresponding sets must be taken.

\begin{corollary}
  Given $(t,x)\in[0,T]\times\RR^n$, the value function $U(t,x)$ is a solution of 
  the equation
  $$
  \sup_{\zeta\in C^+\bs\{0\}}\left[U_{t,\zeta}(t,x)-_\zeta\mathcal{L}_t^*(t,
    -\nabla u_\zeta(t,x),\zeta)-w(t,x,\zeta)\right]=C\, .
  $$
\end{corollary}

\section{Conclusions}

Two features are often present in economic problems. One is that 
the optimization required involves more than one function, which can be in 
conflict.  It is obvious that a linearization is a very strong simplification of the 
problem.  A cone can help to model the preferences of an agent (or the probabilities 
of using one optimization function or the other one).  An approach to the multicriteria case was proposed in \cite{hv} in the classical framework.  The second feature is 
a discount factor that permits to determine the current value of the Lagrangian.  
This was done in \cite{RZ} for the real-valued case.

This paper aims at both targets:  it gives a mathematical model that includes the 
discount factor and the multiobjective nature of the problem.

There are still many open problems.  One of them is the question how to use these 
techniques to handle control problems.  Another line of reasearch could be to 
generalize results like \cite{BJL} to the multicriteria case.

\section*{Acknowledgement}

The author is indebted to Andreas Hamel for fruitful suggestions and
discussions.

This work was supported within the project Optimal Control Problems with 
Set-valued Objective Function by Free University of Bozen-Bolzano (Grant 
OptiConSOF).



\begin{thebibliography}{99}

\bibitem{Aubin} \textsc{J.-P.~Aubin,} \emph{Lax-Hopf formula and Max-Plus 
  properties of solutions to Hamilton-Jacobi equations}, Nonlinear Differential 
  Equations and Applications NoDEA \textbf{20} (2013), 187--211.

\bibitem{Aubin_Frank} \textsc{J.P.~Aubin, H.~Frankowska,} \emph{Set-Valued 
  Analysis,} Birkh\"auser, Boston-Basel-Berlin 1990.

\bibitem{Aumann} \textsc{R.J.~Aumann,} \emph{Integrals of set-valued functions,} 
  J. Math. Anal. Appl.  \textbf{12} (1965), 1--12.

\bibitem{AvantaggiatiLoreti} \textsc{A.~Avantaggiati, P.~Loreti,} 
  \emph{Lax type formulas with lower semicontinuous initial data and 
  hypercontractivity results}, Nonlinear Differential Equations and Applications 
  NoDEA \textbf{20} (2013), 385--411.

\bibitem{Barro} \textsc{R.~Barro,} \emph{Ramsey meets Laibson in the 
  neoclassical growth model}, The Quarterly Journal of Economics \textbf{114} 
  (1999), 1125--1152.

\bibitem{BJL} \textsc{E.N.~Barron,R.~Jensen, W.~Liu,} \emph{Hopf--Lax-Type 
  Formula for $u_t+H(u,Du)=0$}, Journal of Differential Equations \textbf{126} 
  (1996), 48--61.

\bibitem{BV} \textsc{S.P.~Boyd, L.~Vandenberghe,} \emph{Convex 
  Optimization}, Cambridge University Press (2004), 50--51.

\bibitem{Caputo} \textsc{M.~Caputo,} \emph{Foundations of Dynamic 
  Economic Analysis: Optimal Control Theory and Applications}, 
  Cambridge University Press (2005).

\bibitem{DGZ} \textsc{G.~De Marcio, G.~Gorni, G.~Zampieri}, \emph{Global 
  inversion of functions: an introduction}, Nonlinear Differential Equations and 
  Applications \textbf{1} (1994), 229--248.

\bibitem{SetOptSurvey} \textsc{A.H.~Hamel, F.~Heyde, A.~L\"ohne, B.~Rudloff, 
C.~Schrage,} \emph{Set optimization -- A rather short introduction.} In: A.H.~Hamel, F.~Heyde, A.~L\"ohne, B.~Rudloff, C.~Schrage (Eds.), Set 
Optimization and Applications -- The State of the Art, Springer 2015, 65--141.

\bibitem{AC} \textsc{A.H.~Hamel, C.~Schrage,} \emph{Directional derivatives, 
  subdifferentials and optimality conditions for set-valued convex functions,} 
  Pacific Journal of Optimization \textbf{10} 4 (2014), 667--689.

\bibitem{hv} \textsc{A.H.~Hamel, D.~Visetti,} \emph{The value functions approach 
  and Hopf-Lax formula for multiobjective costs via set optimization}, Journal of 
  Mathematical Analysis and Applications \textbf{483} 1 (2020), 123605.

\bibitem{HL} \textsc{F.~Heyde, A.~L\"ohne,} \emph{Solution concepts in vector 
  optimization: a fresh look at an old story}, Optimization \textbf{60} 12 (2011), 
  1421--1440.

\bibitem{Hoang} \textsc{N.~Hoang,} \emph{Hopf-Lax formula and generalized 
  characteristics}, Applicable Analysis \textbf{96} 2 (2013).

\bibitem{Hopf65} \textsc{E.~Hopf,} \emph{Generalized solutions of non-linear 
  equations of first order}, Journal of Mathematics and Mechanics \textbf{14} 
  (1965), 951--973.

\bibitem{Karp} \textsc{L.~Karp,} \emph{Non-constant discounting in continuous 
  time}, Journal of Economic Theory \textbf{132} (2007), 557--568.

\bibitem{Lax57} \textsc{P.D.~Lax,} \emph{Hyperbolic systems of conservation 
  laws II}, Communications on Pure and Applied Mathematics \textbf{10} (1957), 
  537--566.

\bibitem{MSolano_Navas} \textsc{J.~Mar\'\i n-Solano, J.~Navas,} 
  \emph{Non-constant discounting in finite horizon: The free terminal time case}, 
  Journal of Economic Dynamics and Control \textbf{33} (2009), 666--675.

\bibitem{M_SP} \textsc{J.~Mar\'\i n-Solano, C.~Patxot,} \emph{Heterogeneous 
  discounting in economic problems}, Optimal Control Applications and Methods 
  \textbf{33} 1 (2012), 32--50.

\bibitem{M-SS} \textsc{J.~Mar\'\i n-Solano, E.V.~Shevkoplyas,} 
  \emph{Non-constant discounting and differential games with random time 
  horizon}, Automatica \textbf{47} (2011), 2626--2638.

\bibitem{Pilecka} \textsc{M.~Pilecka,} \emph{Set-valued optimization and its 
  application to bilevel optimization,} PhD-thesis 2016, Technische Universit\"at 
  Bergakademie Freiberg.

\bibitem{RZ} \textsc{J.P.~Rinc\'on-Zapatero}, \emph{Hopf-Lax formula for 
  variational problems with non-constant discount}, Journal of Geometric Mechanics   
  \textbf{1} 3 (2009), 357--367.

\bibitem{Stroemberg} \textsc{T.~Str{\"o}mberg,} \emph{The Hopf-Lax formula 
  gives the unique viscosity solution},  Differential Integral Equations \textbf{15} 1 
  (2002), 47--52.

\end{thebibliography}
\end{document}